\documentclass[11pt, reqno, psamsfonts]{amsart}
\pdfoutput=1

\usepackage{amssymb}
\usepackage{amsthm}
\usepackage{amsmath}
\usepackage{latexsym}
\usepackage[T2A,OT1]{fontenc}
\usepackage[utf8]{inputenc}
\usepackage[russian, french, english]{babel}

\usepackage{graphicx}
\usepackage{wrapfig}
\usepackage[justification=centering, labelfont=bf]{caption}
\usepackage{subcaption}
\usepackage{mathtools}
\usepackage[hidelinks]{hyperref}
\usepackage{enumitem}
\usepackage{mathrsfs}
\usepackage{float}
\usepackage{stmaryrd}
\usepackage{framed}
\usepackage{tikz}
\usepackage{lmodern}
\usepackage{mathabx}
\usepackage[shortcuts]{extdash}
\usepackage[foot]{amsaddr}
\usepackage{titlesec}
\usetikzlibrary{shapes,snakes}
\usetikzlibrary{arrows.meta}

\usepackage[spacing=true,kerning=true,babel=true,tracking=true]{microtype}

\usepackage[backend=biber, style=alphabetic, sorting=nyt, maxnames=100]{biblatex}
\usepackage[left=2.3cm,right=2.3cm,top=2.3cm,bottom=2.3cm,bindingoffset=0cm]{geometry}

\title{Fractional DP-colorings of Sparse Graphs}
\date{}
\author{Anton~Bernshteyn}
\address[Anton Bernshteyn]{Department of Mathematics, University of Illinois at Urbana--Champaign, IL, USA and Department of Mathematical Sciences, Carnegie Mellon University, Pittsburgh, PA, USA}
\email{abernsht@math.cmu.edu}

\author{Alexandr Kostochka}
\address[Alexandr Kostochka]{Department of Mathematics, University of Illinois at Urbana--Champaign, IL, USA and Sobolev Institute of Mathematics, Novosibirsk, Russia}
\email{kostochk@math.uiuc.edu}
\thanks{Research of the second author is supported in part by NSF grant
	DMS-1600592 and grants 18-01-00353  and 19-01-00682 of the Russian Foundation for Basic Research.}

\author{Xuding Zhu}
\address[Xuding Zhu]{Department of Mathematics, Zhejiang Normal University, Jinhua, China}
\thanks{Research of the third author is supported in part by CNSF grant 11571319.}
\email{xudingzhu@gmail.com}

\newtheoremstyle{bfnote}%
{}{}%
{\slshape}{}%
{\bfseries}{\bfseries.}%
{ }%
{\thmname{#1}\thmnumber{ #2}\thmnote{ \ep{\normalfont{}#3}}}

\newtheoremstyle{defbfnote}%
{}{}%
{}{}%
{\bfseries}{.}%
{ }%
{\thmname{#1}\thmnumber{ #2}\thmnote{ (#3)}}

\newtheoremstyle{claim}%
{}{}%
{\slshape}{}%
{\itshape}{.}%
{ }%
{\thmname{#1}\thmnumber{ #2}\thmnote{ \ep{\normalfont{}#3}}}

\theoremstyle{bfnote}
\newtheorem{theo}[equation]{Theorem}

\newtheorem{lemma}[equation]{Lemma}
\newtheorem{corl}[equation]{Corollary}

\newtheorem*{obs}{Observation}
\newtheorem*{claim*}{Claim}
\newtheorem*{corl*}{Corollary}

\newcommand*{\myproofname}{Proof}
\newenvironment{claimproof}[1][\myproofname]{\begin{proof}[#1]}{\end{proof}}

\theoremstyle{definition}
\newtheorem{defn}[equation]{Definition}
\newtheorem*{defn*}{Definition}

\newtheorem{ques}[equation]{Question}

\newtheorem*{exmp*}{Example}

\theoremstyle{remark}
\newtheorem*{ques*}{Question}
\newtheorem*{remk*}{Remark}

\theoremstyle{claim}

\newcounter{ForClaims}[section]
\newtheorem{claim}{Claim}[ForClaims]

\newcommand{\0}{\varnothing}
\newcommand{\set}[1]{\{#1\}}

\newcommand{\N}{{\mathbb{N}}}
\newcommand{\Z}{\mathbb{Z}}
\newcommand{\F}{\mathbb{F}}
\newcommand{\R}{\mathbb{R}}
\renewcommand{\C}{\mathbb{C}}
\newcommand{\Q}{\mathbb{Q}}

\renewcommand{\epsilon}{\varepsilon}
\renewcommand{\phi}{\varphi}
\renewcommand{\theta}{\vartheta}

\renewcommand{\leq}{\leqslant}
\renewcommand{\geq}{\geqslant}

\newcommand{\powerset}[1]{\operatorname{Pow}(#1)}

\renewcommand{\G}{\Gamma}

\newcommand{\defeq}{\coloneqq}

\newcommand{\T}{\mathbb{T}}

\newcommand{\Cov}[1]{\mathscr{#1}}
\newcommand{\bemph}[1]{{\normalfont#1}} 
\newcommand{\ep}[1]{\bemph{(}#1\bemph{)}} 

\newenvironment{scproof}[1][Proof]{\begin{proof}[\textsc{#1}]}{\end{proof}}

\numberwithin{equation}{section}

\renewcommand{\thesubsection}{\arabic{section}.\Alph{subsection}}

\usepackage{etoolbox}

\titleformat{\section}[block]{\scshape\filcenter}{\thesection.}{1ex}{}
\titleformat{\subsection}[block]{\bfseries\filcenter}{\thesubsection.}{1ex}{}
\titleformat{\subsubsection}[runin]{\bfseries}{\thesubsubsection.}{1ex}{}[.]

\titlespacing*{\section}{0pt}{*3}{*1}
\titlespacing*{\subsection}{0pt}{*2}{*1}

\makeatletter
\newcommand{\neutralize}[1]{\expandafter\let\csname c@#1\endcsname\count@}
\makeatother

\renewcommand{\mkbibnamefirst}[1]{\textsc{#1}}
\renewcommand{\mkbibnamelast}[1]{\textsc{#1}}
\renewcommand{\mkbibnameprefix}[1]{\textsc{#1}}
\renewcommand{\mkbibnameaffix}[1]{\textsc{#1}}
\renewcommand{\finentrypunct}{}

\renewbibmacro{in:}{}

\renewbibmacro*{volume+number+eid}{%
	\printfield{volume}%
	\setunit*{\addnbspace}
	\printfield{number}%
	\setunit{\addcomma\space}%
	\printfield{eid}}

\DeclareFieldFormat[article]{volume}{\textbf{#1}\space}
\DeclareFieldFormat[article]{number}{\mkbibparens{#1}}

\DeclareFieldFormat{journaltitle}{#1,}
\DeclareFieldFormat[thesis]{title}{\mkbibemph{#1}\addperiod}
\DeclareFieldFormat[article, unpublished, thesis]{title}{\mkbibemph{#1},}
\DeclareFieldFormat[book]{title}{\mkbibemph{#1}\addperiod}
\DeclareFieldFormat[unpublished]{howpublished}{#1, }

\DeclareFieldFormat{pages}{#1}

\DeclareFieldFormat[article]{series}{Ser.~#1\addcomma}

\begin{document}
	\pagestyle{plain}
	
	\maketitle
	
	\begin{abstract}
		DP-coloring (also known as correspondence coloring) is a generalization of list coloring developed recently by Dvo\v r\' ak and Postle~\cite{DP15}. In this paper we introduce and study the fractional DP-chromatic number $\chi_{DP}^\ast(G)$. We characterize all connected graphs $G$ such that $\chi_{DP}^\ast(G) \leq 2$: they are precisely the graphs with no odd cycles and at most one even cycle. By a theorem of Alon, Tuza, and Voigt~\cite{ATV1997}, the fractional list-chromatic number $\chi_\ell^\ast(G)$ of any graph $G$ equals
		its fractional chromatic number $\chi^\ast(G)$. This equality does not extend to fractional DP-colorings. Moreover, we show that the difference $\chi^\ast_{DP}(G) - \chi^\ast(G)$ can be arbitrarily large, and, furthermore, $\chi^\ast_{DP}(G) \geq d/(2 \ln d)$ for every graph $G$ of maximum average degree $d \geq 4$. On the other hand, we show that this asymptotic lower bound is tight for a large class of graphs that includes all bipartite graphs as well as many graphs of high girth and high chromatic number.
	\end{abstract}
	
	\section{Introduction}
	
	\noindent 	All graphs considered here are finite and do not have multiple edges and loops. By a ``graph'' we mean an undirected graph; directed graphs are referred to as \emph{digraphs}. 
	
	 \emph{DP-coloring}, also known as \emph{correspondence coloring}, is a generalization of list coloring introduced recently by Dvo\v r\' ak and Postle~\cite{DP15}. In the setting of DP-coloring, not only does each vertex get its own list of available colors, but also the identifications between the colors in the lists are allowed to vary from edge to edge.
	
	\begin{defn}\label{defn:cover}
		Let $G$ be a graph. A \emph{cover} of $G$ is a pair $\Cov{H} = (L, H)$, consisting of a graph $H$ and a function $L \colon V(G) \to \powerset{V(H)}$, satisfying the following requirements:
		\begin{enumerate}[labelindent=\parindent,leftmargin=*,label=(C\arabic*)]
			\item the sets $\set{L(u) \,:\,u \in V(G)}$ form a partition of $V(H)$;
			\item for every $u \in V(G)$, the graph $H[L(u)]$ is complete;
			\item if $E_H(L(u), L(v)) \neq \0$, then either $u = v$ or $uv \in E(G)$;
			\item \label{item:matching} if $uv \in E(G)$, then $E_H(L(u), L(v))$ is a matching.
		\end{enumerate}
		A cover $\Cov{H} = (L, H)$ of $G$ is \emph{$k$-fold} if $|L(u)| = k$ for all $u \in V(G)$.
	\end{defn}
	
	\begin{remk*}
		We emphasize that the matching $E_H(L(u), L(v))$ in Definition~\ref{defn:cover}\ref{item:matching} is not required to be perfect and, in particular, may be empty.
	\end{remk*}
	
	
	
	\begin{defn}
		Let $G$ be a graph and let $\Cov{H} = (L, H)$ be a cover of $G$. An \emph{$\Cov{H}$-coloring} of $G$ is an independent set in $H$ of size $|V(G)|$. The \emph{DP-chromatic number} $\chi_{DP}(G)$ of $G$ is the smallest $k \in \N$ such that $G$ admits an $\Cov{H}$-coloring for every $k$-fold cover $\Cov{H}$ of $G$.
	\end{defn}
	
	\begin{remk*}\label{remk:single}
		By definition, if $\Cov{H} = (L, H)$ is a cover of a graph $G$, then $\set{L(u)\,:\, u \in V(G)}$ is a partition of~$H$ into $|V(G)|$ cliques. Therefore, an independent set $I \subseteq V(H)$ is an $\Cov{H}$-coloring of $G$ if and only if $|I \cap L(u)| = 1$ for all $u \in V(G)$.
	\end{remk*}

		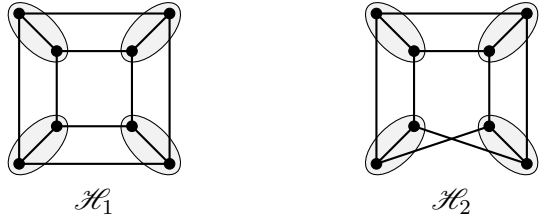
\begin{figure}[t!]
		\centering	
		\begin{tikzpicture}[scale=0.5]
		\definecolor{light-gray}{gray}{0.95}
		
		\filldraw[fill=light-gray]
		(6.5,0) circle [x radius=1cm, y radius=5mm, rotate=45]
		(6.5,3) circle [x radius=1cm, y radius=5mm, rotate=-45]
		(9.5,0) circle [x radius=1cm, y radius=5mm, rotate=-45]
		(9.5,3) circle [x radius=1cm, y radius=5mm, rotate=45];
		
		\foreach \x in {(6, -0.5), (6, 3.5), (7, 0.5), (7, 2.5), (9, 0.5), (9, 2.5), (10, -0.5), (10, 3.5)}
		\filldraw \x circle (4pt);
		
		\draw[thick] (6, -0.5) -- (6, 3.5) -- (10, 3.5) -- (10, -0.5) -- cycle;
		
		\draw[thick] (7, 0.5) -- (7, 2.5) -- (9, 2.5) -- (9, 0.5) -- cycle;
		
		\draw[thick] (6, -0.5) -- (7, 0.5) (6, 3.5) -- (7, 2.5) (9, 2.5) -- (10, 3.5) (10, -0.5) -- (9, 0.5);
		
		\node at (8, -1.5) {$\Cov{H}_1$};
		
		\filldraw[fill=light-gray]
		(13+3,0) circle [x radius=1cm, y radius=5mm, rotate=45]
		(13+3,3) circle [x radius=1cm, y radius=5mm, rotate=-45]
		(16+3,0) circle [x radius=1cm, y radius=5mm, rotate=-45]
		(16+3,3) circle [x radius=1cm, y radius=5mm, rotate=45];
		
		\foreach \x in {(12.5+3, -0.5), (12.5+3, 3.5), (13.5+3, 0.5), (13.5+3, 2.5), (15.5+3, 0.5), (15.5+3, 2.5), (16.5+3, -0.5), (16.5+3, 3.5)}
		\filldraw \x circle (4pt);
		
		\draw[thick] (12.5+3, -0.5) -- (12.5+3, 3.5) -- (16.5+3, 3.5) -- (16.5+3, -0.5) -- (13.5+3, 0.5) -- (13.5+3, 2.5) -- (15.5+3, 2.5) -- (15.5+3, 0.5) -- cycle;
		
		\draw[thick] (12.5+3, -0.5) -- (13.5+3, 0.5) (12.5+3, 3.5) -- (13.5+3, 2.5) (15.5+3, 0.5) -- (16.5+3, -0.5) (15.5+3, 2.5) -- (16.5+3, 3.5);
		
		\node at (14.5+3, -1.5) {$\Cov{H}_2$};
		\end{tikzpicture}
		\caption{Two distinct $2$-fold covers of the $4$-cycle.
		}\label{fig:cycle}
	\end{figure}
	
	
		As an illustration, consider the two $2$-fold covers of the $4$-cycle $C_4$ shown in Fig.~\ref{fig:cycle}. Observe that $C_4$ is $\Cov{H}_1$-colorable but not $\Cov{H}_2$-colorable; in particular, $\chi_{DP}(C_4) \geq 3$. On the other hand, it is easy to see that $\chi_{DP}(G) \leq \Delta + 1$ for any graph $G$ of maximum degree $\Delta$, so $\chi_{DP}(C_4) = 3$. The same argument shows that $\chi_{DP}(C_n) = 3$ for every cycle $C_n$. 
	
	To see that DP-coloring is a generalization of list coloring, suppose that $G$ is a graph and $L$ is a list assignment for $G$. Let $H$ be the graph with vertex set
	\[
	V(H) \defeq \set{(u, c)\,:\, u \in V(G) \text{ and } c \in L(u)},
	\]
	in which two distinct vertices $(u, c)$ and $(v, d)$ are adjacent if and only if
	\begin{itemize}
		\item[--] either $u = v$,
		\item[--] or else, $uv \in E(G)$ and $c = d$.
	\end{itemize}
	For each $u \in V(G)$, set
	$
		L'(u) \defeq \set{(u, c) \,:\, c \in L(u)}
	$.
	Then $\Cov{H} \defeq (L', H)$ is a cover of $G$, and there is a natural bijective correspondence between the $L$-colorings and the $\Cov{H}$-colorings of $G$. This, in particular, implies that $\chi_{DP}(G) \geq \chi_\ell(G)$ for all $G$.
	
	In this paper we introduce and study the fractional version of DP-coloring. We start with a brief review of the classical concepts of fractional coloring and fractional list coloring. For a survey of the topic, see, e.g., \cite[Chapter~3]{FracBook}.
	
	Let $G$ be a graph. An \emph{$(\eta, k)$-coloring} of $G$, where $\eta \in [0,1]$ and $k \in \N^+$, is a map $f \colon V(G) \to \powerset{[k]}$ with the following properties:
	\begin{enumerate}[labelindent=\parindent,leftmargin=*,label=(F\arabic*)]
		\item for every vertex $u \in V(G)$, we have $|f(u)| \geq \eta k$;
		\item for every edge $uv\in E(G)$, we have $f(u) \cap f(v) = \0$.
	\end{enumerate}
	For given $k \in \N^+$, let
	\[
		\theta(G, k) \defeq \max \set{\eta \in [0,1] \,:\, \text{$G$ admits an $(\eta, k)$-coloring}}.
	\]
	(The maximum is attained, as only the values of the form $\ell/k$ for integer $\ell$ are relevant.) The \emph{fractional chromatic number} $\chi^\ast(G)$ of $G$ is defined by
	\begin{equation}\label{eq:frac}
	\chi^\ast(G) \defeq \inf \set{\theta(G, k)^{-1} \,:\, k \in \N^+}.
	\end{equation}
	It is well-known~\cite[\S3.1]{FracBook} that the infimum in \eqref{eq:frac} is actually a minimum: For every graph $G$, there is some $k \in \N^+$ such that $\chi^\ast(G) = \theta(G, k)^{-1}$. In particular, $\chi^\ast(G)$ is always a rational number.
	
	Fractional coloring allows a natural list-version. Let $G$ be a graph and let $L$ be a list assignment for $G$. An \emph{$(\eta, L)$-coloring} of $G$, where $\eta \in [0,1]$, is a map $f$ that associates to each $u \in V(G)$ a subset $f(u) \subseteq L(u)$ with the following properties:
	\begin{enumerate}[labelindent=\parindent,leftmargin=*,label=(FL\arabic*)]
		\item for every vertex $u \in V(G)$, we have $|f(u)| \geq \eta |L(u)|$;
		\item for every edge $uv\in E(G)$, we have $f(u) \cap f(v) = \0$.
	\end{enumerate}
	We say that $L$ is a \emph{$k$-list assignment} if $|L(u)| = k$ for all $u \in V(G)$. For given $k \in \N^+$, let
	\[
		\theta_\ell(G, k) \defeq \max \set{\eta \in [0,1] \,:\, \text{$G$ admits an $(\eta, L)$-coloring for every $k$-list assignment $L$ for $G$}}.
	\]
	The \emph{fractional list-chromatic number} $\chi_\ell^\ast(G)$ of $G$ is defined by
	\[
		\chi_\ell^\ast(G) \defeq \inf \set{\theta_\ell(G, k)^{-1} \,:\, k \in \N^+}.
	\]
	Somewhat surprisingly, Alon, Tuza, and Voigt~\cite{ATV1997} showed that $\chi_\ell^\ast(G) = \chi^\ast(G)$ for all graphs $G$ and, in fact, for each $G$, there is $k \in \N^+$ such that
	\[
		\chi^\ast_\ell(G) = \chi^\ast(G) =  \theta_\ell(G, k)^{-1} = \theta(G, k)^{-1}.
	\]
	 (Recall that the list-chromatic number of a graph cannot be bounded above by any function of its ordinary chromatic number.)
	
	Now we proceed to our main definitions. Given a cover $\Cov{H} = (L, H)$ of a graph $G$, we refer to the edges of $H$ connecting distinct parts of the partition $\set{L(u) \,:\, u \in V(G)}$ as \emph{cross-edges}. A subset $S \subseteq V(H)$ is \emph{quasi-independent} if it spans no cross-edges. 
	
	\begin{defn}
		Let $\Cov{H} = (L, H)$ be a cover of a graph $G$ and let $\eta \in [0,1]$. An \emph{$(\eta, \Cov{H})$-coloring} of $G$ is a quasi\=/independent set $S \subseteq V(H)$ such that $|S \cap L(u)| \geq \eta |L(u)|$ for all $u \in V(G)$. 
	\end{defn}
	
	\begin{defn}
		Let $G$ be a graph. For $k \in \N^+$, let
		\[
			\theta_{DP}(G, k) \defeq \max \set{\eta \in [0,1] \,:\, \text{$G$ admits an $(\eta, \Cov{H})$-coloring for every $k$-fold cover $\Cov{H}$ of $G$}}.
		\]
		The \emph{fractional DP-chromatic number} $\chi_{DP}^\ast(G)$ is defined by
		\begin{equation}\label{eq:frac_DP}
			\chi_{DP}^\ast(G) \defeq \inf \set{\theta_{DP}(G,k)^{-1} \,:\, k \in \N^+}.
		\end{equation}
	\end{defn}
	
	Clearly, for each graph $G$, we have $\chi^\ast(G) \leq \chi_{DP}^\ast(G) \leq \chi_{DP}(G)$. Our results described below imply that both inequalities can be strict.
	
	Since $\chi_{DP}(C_n) = 3$ for every cycle $C_n$, a connected graph $G$ satisfies $\chi_{DP}(G) \leq 2$ if and only if $G$ is a tree. Our first result is the characterization of graphs $G$ with $\chi_{DP}^\ast(G) \leq 2$:
	
	\begin{theo}\label{theo:2}
		Let $G$ be a connected graph. Then $\chi_{DP}^\ast(G) \leq 2$ if and only if $G$ contains no odd cycles and at most one even cycle. Furthermore, if $G$ contains no odd cycles and exactly one even cycle, then \[\chi_{DP}^\ast(G) = 2, \quad \text{even though} \quad \theta_{DP}(G, k)^{-1} > 2 \text{ for all } k \in \N^+\] \ep{i.e., the infimum in \eqref{eq:frac_DP} is not attained}.
	\end{theo}
	
	Theorem~\ref{theo:2} shows that the Alon--Tuza--Voigt theorem does not extend to fractional DP-coloring, as every connected bipartite graph $G$ with $|E(G)| \geq |V(G)| + 1$ satisfies $\chi^\ast(G) = \chi(G) = 2$, while $\chi_{DP}^\ast(G) > 2$. Theorem~\ref{theo:2} also provides examples of graphs for which the infimum in \eqref{eq:frac_DP} is not attained. However, the following natural question remains open:
	
	\begin{ques}
		Do there exist graphs $G$ for which $\chi_{DP}^\ast(G)$ is irrational?
	\end{ques}
	
	The \emph{average degree} of a graph $G$ with $V(G) \neq \0$ is the value $2|E(G)|/|V(G)|$; the \emph{maximum average degree} of $G$ is the largest average degree of a subgraph $H$ of $G$ with $V(H) \neq \0$. It was shown in~\cite{Ber} that $\chi_{DP}(G) = \Omega(d/\ln d)$ for graphs $G$ of maximum average degree $d$. Using a similar argument, we extend this asymptotic lower bound to the fractional setting:
	
	\begin{theo}\label{theo:lower_bound}
		If $G$ is a graph of maximum average degree $d \geq 4$, then $\chi^\ast_{DP}(G) \geq d/(2\ln d)$.
	\end{theo}
	
	From Theorem~\ref{theo:lower_bound}, it follows that $\chi^\ast_{DP}(G)$ cannot be bounded above by any function of $\chi^\ast(G)$, since there exist bipartite graphs of arbitrarily high average degree.
	
	
	Recall that a graph $G$ is \emph{$d$-degenerate} if it has an orientation $D$ such that:
	\begin{itemize}
		\item[--] $D$ is \emph{acyclic}, i.e., there are no directed cycles in $D$;
		
		\item[--] $\Delta^+(D) \leq d$, where $\Delta^+(D)$ denotes the maximum out-degree of $D$.
	\end{itemize}
	Note that every graph of maximum average degree $d$ is $d$-degenerate, while every $d$-degenerate graph has maximum average degree at most $2d$. Our next result describes additional conditions on an acyclic orientation $D$ under which the lower bound given by Theorem~\ref{theo:lower_bound} is asymptotically tight. Throughout, given a digraph $D$, we write $E(D)$ for the set of all \emph{directed} edges of $D$ \ep{so $uv \in E(D)$ indicates a directed edge from $u$ to $v$}.
	
	\begin{theo}\label{theo:frac_deg}
		Suppose that a graph $G$ has an acyclic orientation $D$ such that
		\begin{enumerate}[label=\normalfont{(D\arabic*)}]
			\item $\Delta^+(D) \leq d$; and
			\item\label{item:parity} for all $uv \in E(D)$, there is no directed $uv$-path of even length in $D$.
		\end{enumerate}
		Then $\chi_{DP}^\ast(G) \leq (1 + o(1))d/\ln d$.
	\end{theo}
	
	Obviously, every orientation $D$ of a \emph{bipartite} graph $G$ satisfies condition~\ref{item:parity} of Theorem~\ref{theo:frac_deg}. Hence, we obtain the following:
	
	\begin{corl}\label{corl:bip}
		If $G$ is a $d$-degenerate bipartite graph, then $\chi_{DP}^\ast(G) \leq (1 + o(1))d/\ln d$.
	\end{corl}

	Corollary~\ref{corl:bip} shows that the lower bound given by Theorem~\ref{theo:lower_bound} is sharp, including the value of the constant factor. For example, if $\omega \colon \N \to \N$ is any function such that $\omega(d)/d \to \infty$ as $d \to \infty$, then the complete bipartite graph $K_{d,\,\omega(d)}$ is $d$-degenerate, while its average degree is $(2+o(1))d$. Thus, Theorem~\ref{theo:lower_bound} and Corollary~\ref{corl:bip} combine to give
	$
		\chi^\ast_{DP} (K_{d,\, \omega(d)}) = (1+o(1)) d/\ln d
	$.
	
	The conclusion of Theorem~\ref{theo:frac_deg} is interesting even for the ordinary fractional chromatic number, especially since its requirements are satisfied by several known constructions of graphs with high girth and high chromatic number. One example is the scheme analyzed in~\cite{KN99}, based on the Blanche Descartes construction \cite{BD} of triangle-free graphs with high chromatic number. For completeness, we sketch it in Section~\ref{sec:BD}. Another related family of graphs of high chromatic number that falls under the conditions of Theorem~\ref{theo:frac_deg} is described in~\cite[Theorem 3.4]{AKRWZ}. From these constructions, we deduce the following result:
	
	\begin{corl}\label{corl:large_girth}
		For all $d$, $g \in \N^+$, there exists a graph $G_{d, g}$ with chromatic number at least $d$, girth at least $g$, and $\chi^\ast_{DP}(G_{d,g}) \leq (1 + o(1))d/\ln d$.
	\end{corl}

	The remainder of this paper is organized as follows. First, we prove Theorem~\ref{theo:lower_bound} in Section~\ref{sec:lower_bound}. Then, in Section~\ref{sec:2}, we establish Theorem~\ref{theo:2}. Section~\ref{sec:frac_deg} is dedicated to the proof of Theorem~\ref{theo:frac_deg}. Finally, the proof of Corollary~\ref{corl:large_girth} is given in Section~\ref{sec:BD}.
	
	\section{Proof of Theorem~\ref{theo:lower_bound}}\label{sec:lower_bound}
	
	\noindent What follows is a slight modification of~\cite[Theorem~1.6]{Ber}. Let $G$ be a graph of maximum average degree $d \geq 4$. After passing to a subgraph, we may assume that the average degree of $G$ itself is $d$. Set $n \defeq |V(G)|$ and $m \defeq |E(G)|$ \ep{and thus $m = dn/2$}. Let $\eta_0 \defeq 2\ln d/ d$. Our goal is to show that $\theta_{DP}(G, k) < \eta_0$ for all $k \in \N^+$.
	To that end, fix arbitrary $k \in \N^+$ and let $\eta \defeq \lceil \eta_0 k \rceil / k$. It is enough to prove $\theta_{DP}(G, k) < \eta$.
	
	Let $\set{L(u)\,:\, u \in V(G)}$ be a collection of pairwise disjoint $k$-element sets. Define $X \defeq \bigcup_{u \in V(G)} L(u)$ and build a random graph $H$ with vertex set $X$ by making each $L(u)$ a clique and putting, independently for each $uv \in E(G)$, a uniformly random perfect matching between $L(u)$ and $L(v)$. Let $\Cov{H} \defeq (L, H)$ denote the resulting random $k$-fold cover of $G$. We shall argue that, with positive probability, $G$ is not $(\eta, \Cov{H})$-colorable.
	
	Let $S \subseteq X$ be an arbitrary set with $|S \cap L(u)| = \eta k$ for all $u \in V(G)$. Consider any edge $uv \in E(G)$. The $H$-neighborhood of the set $S \cap L(u)$ in $L(v)$ is a uniformly random $(\eta k)$-element subset of $L(v)$, so it is disjoint from $S \cap L (v)$ with probability
	\[
		{(1-\eta)k \choose \eta k }{k \choose \eta k}^{-1}.
	\]
	Since the matchings corresponding to different edges of $G$ are drawn independently from each other, we conclude that
	\begin{align*}
	&\Pr \left[S \text{ is quasi-independent in } H\right] \\ = & \prod_{uv \,\in\, E(G)} \Pr\left[\text{there are no cross-edges between $S \cap L(u)$ and $S \cap L(v)$}\right]
	= {(1-\eta)k \choose \eta k }^m{k \choose \eta k}^{-m}.
	\end{align*}
	There are ${k \choose \eta k}^n$ possible choices for $S$, so we can use the union bound to get
	\[
	\Pr\left[\text{$G$ is $(\eta, \Cov{H})$-colorable}\right] \leq {(1-\eta)k \choose \eta k }^m{k \choose \eta k}^{-m} {k \choose \eta k}^{n} = \left({(1-\eta)k \choose \eta k}^{d/2} {k \choose \eta k}^{-(d/2-1)} \right)^{n}.
	\]
	Thus, we only need to show that
	\[
	{(1-\eta)k \choose \eta k}^{d/2} {k \choose \eta k}^{-(d/2-1)} < 1.
	\]
	Notice that
	\[
	{(1-\eta)k \choose \eta k} {k \choose \eta k}^{-1} = \prod_{i=0}^{\eta k -1} \frac{(1-\eta)k - i}{k - i} \leq (1-\eta)^{\eta k}.
	\]
	Additionally,
	\[
	{k \choose \eta k} \leq \left(\frac{e}{\eta}\right)^{\eta k}.
	\]
	Therefore,
	\[
	{(1-\eta)k \choose \eta k}^{d/2} {k \choose \eta k}^{-(d/2-1)} \leq \left(\frac{e(1-\eta)^{d/2}}{\eta}\right)^{\eta k},
	\]
	so it is enough to establish
	\[
	e(1-\eta)^{d/2} < \eta.
	\]
	Since $1-\eta \leq \exp(-\eta)$, we have
	\[
	e(1-\eta)^{d/2} \leq e \cdot \exp(-\eta d/2) \leq e d^{-1} < \eta,
	\]
	as long as $d > e^{e/2} \approx 3.89$, as desired.
	
	\section{Proof of Theorem~\ref{theo:2}}\label{sec:2}
	
	\begin{lemma}\label{lemma:nplus1}
		If $G$ is a graph such that $|E(G)| \geq |V(G)| + 1$, then $\chi_{DP}^\ast(G) > 2$.
	\end{lemma}
	\begin{scproof}
		Set $n \defeq |V(G)|$. Without loss of generality, we may assume that $|E(G)| = n + 1$. Pick a positive real number $\eta_0 < 1/2$ so that for all $\eta \in [\eta_0, 1/2)$, the following inequality holds:
		\[
			\left(\frac{e(1-\eta)}{1-2\eta}\right)^{(1-2\eta) (n+1)} < \left(\frac{1}{\eta}\right)^{\eta}.
		\]
		Such $\eta_0$ exists since, as $\eta$ approaches $1/2$ from below, we have
		\[
		\left(\frac{e(1-\eta)}{1-2\eta}\right)^{(1-2\eta)(n+1)} \to 1, \qquad \text{while} \qquad \left(\frac{1}{\eta}\right)^{\eta} \to \sqrt{2}.
		\] 
		Our aim is to show that for all $k \in \N^+$, $\theta_{DP}(G, k) < \eta_0$. Fix $k \in \N^+$ and let $\eta \defeq \lceil\eta_0 k\rceil/k$, so it suffices to show that $\theta_{DP}(G, k) < \eta$. We may assume that $\eta \leq 1/2$, as $\theta_{DP}(G, k) \leq 1/2$ for every graph $G$ with at least one edge. We use the same approach and notation as in the proof of Theorem~\ref{theo:lower_bound} \ep{see Section~\ref{sec:lower_bound}}. Thus, $\Cov{H} = (L, H)$ is a random $k$-fold cover of $G$, where $V(H) = X$, and if $S \subseteq X$ is a set with $|S \cap L(u)| = \eta k$ for all $u \in V(G)$, then
		\[
		\Pr \left[S \text{ is quasi-independent in } H\right] = {(1-\eta)k \choose \eta k}^{n+1}{k \choose \eta k}^{-(n+1)},
		\]
		so the probability that $G$ is $(\eta, \Cov{H})$-colorable is at most
		\[
		{(1-\eta)k \choose \eta k}^{n+1}{k \choose \eta k}^{-(n+1)} {k \choose \eta k}^{n} = {(1-\eta)k \choose \eta k}^{n+1} {k \choose \eta k}^{-1}.
		\]
		If $\eta = 1/2$, then the last expression is equal to
		\[
			{k \choose k/2}^{-1} < 1,
		\]
		as desired. Thus, we may assume from now on that $\eta < 1/2$. Then we can write
		\[
			{(1-\eta)k \choose \eta k} = {(1-\eta)k \choose (1-2\eta) k} \leq \left(\frac{e(1-\eta)}{1-2\eta}\right)^{(1-2\eta) k}.
		\]
		Additionally,
		\[
			{k \choose \eta k} \geq \left(\frac{1}{\eta}\right)^{\eta k}.
		\]
		Therefore, the probability that $G$ is $(\eta, \Cov{H})$-colorable is less than $1$ as long as
		\[
			\left(\frac{e(1-\eta)}{1-2\eta}\right)^{(1-2\eta) (n+1)} < \left(\frac{1}{\eta}\right)^{\eta},
		\]
		which is true by the choice of $\eta_0$ and since $\eta_0 \leq \eta < 1/2$.
	\end{scproof}
	
	\begin{lemma}\label{lemma:even}
		If $G$ is a cycle of even length, then $\chi_{DP}^\ast(G) = 2$, while $\theta_{DP}(G, k)^{-1} > 2$ for all $k \in \N^+$.
	\end{lemma}
	\begin{scproof}
		Let the vertex and the edge sets of $G$ be $\set{v_1, \ldots, v_n}$ and $\set{v_1v_2, v_2v_3, \ldots, v_n v_1}$ respectively. Given $k \in \N^+$ and a permutation $\sigma \colon [k] \to [k]$, we define a $k$-fold cover $\Cov{H}_\sigma = (L_\sigma, H_\sigma)$ of $G$ as follows. First, for each $1 \leq i \leq k$, let $L_\sigma(v_i) \defeq \set{i} \times [k]$. Then, for each $1 \leq i < n$, define
		\[
			E_{H_\sigma}(L_\sigma(v_i), L_\sigma(v_{i+1})) \defeq \set{\set{(i, j), (i+1, j)} \,:\, 1 \leq j \leq k}.
		\]
		Finally, let
		\[
		E_{H_\sigma} (L_\sigma(v_1), L_\sigma(v_n)) \defeq \set{\set{(1,j), (n, \sigma(j))} \,:\, 1 \leq j \leq k}.
		\]
		It is clear that to determine $\theta_{DP}(G, k)$ it is enough to consider covers of the form $\Cov{H}_\sigma$ for some $\sigma$.
		
		Suppose that $\theta_{DP}(G, k) = 1/2$ for some $k \in \N$. Consider a permutation $\sigma \colon [k] \to [k]$ that consists of a single cycle. Note that if $X \subseteq [k]$ satisfies $\sigma(X) = X$, then $X \in \set{\0, [k]}$. Let $S$ be a $(1/2, \Cov{H}_\sigma)$-coloring of $G$. For each $1 \leq i \leq k$, let $S_i \defeq \set{j \,:\, (i,j) \in S}$. Since $S$ is quasi-independent, $S_i \cap S_{i + 1} = \0$ for all $1 \leq i < n$. But we also have $|S_i| = |S_{i+1}| = k/2$, and thus $S_{i+1} = [k] \setminus S_i$. Since $n$ is even, this yields $S_n = [k] \setminus S_1$. For every $j \in S_1$, we have $\sigma(j) \not \in S_n$, hence $\sigma(j) \in S_1$. In other words, $\sigma(S_1) = S_1$. But then $S_1 \in \set{\0, [k]}$, contradicting the fact that $|S_1| = k/2$.
		
		It remains to prove that for any $\eta < 1/2$, there is $k \in \N^+$ such that $\theta_{DP}(G, k) \geq \eta$. Take a large odd integer $k$ and let $\sigma \colon [k] \to [k]$ be a permutation. Write $\sigma$ as a product of disjoint cycles: $\sigma = \pi_1 \cdots \pi_m$. We may rearrange the set $[k]$ so that the support of each cycle $\pi_i$ is of the form $\set{m \in [k] \,:\, \ell_i \leq m \leq r_i}$ for some $1 \leq \ell_i \leq r_i \leq k$, and
		\[
			\pi_i (\ell_i) = \ell_i + 1, \quad \pi_i (\ell_i + 1) = \ell_i + 2, \quad \ldots, \quad \pi_i(r_i) = \ell_i.
		\]
		Then $\sigma(i) \leq i + 1$ for all $1 \leq i \leq k$. Now let
		\[
			X \defeq \set{1, \ldots, (k-1)/2} \quad \text{and} \quad Y \defeq \set{(k+3)/2, \ldots, k}.
		\]
		Note that $|X|= |Y| = (k-1)/2$, $X \cap Y = \0$, and $\sigma(X) \cap Y = \0$. Hence, if we define
		\[
			S \defeq \set{(i, j) \,:\, 1 \leq i \leq n,\, j \in X \text{ if } i \text{ is odd and } j \in Y \text{ if } i \text{ is even}},
		\]
		then $S$ is a $((1-1/k)/2, \Cov{H}_\sigma)$-coloring of $G$, and we are done.
	\end{scproof}
	
	\begin{scproof}[Proof of Theorem~\ref{theo:2}]
		Let $G$ be a connected graph and suppose that $\chi_{DP}^\ast(G) \leq 2$. Even the ordinary fractional chromatic number of any odd cycle exceeds $2$ (see \cite[Proposition~3.1.2]{FracBook}), so $G$ must be bipartite. Furthermore, by Lemma~\ref{lemma:nplus1}, $|E(G)| \leq |V(G)|$, so $G$ contains at most one even cycle. Conversely, suppose that $G$ contains no odd cycles and at most one even cycle. If $G$ is acyclic, then $\chi_{DP}^{\ast}(G) = \chi_{DP}(G) \leq 2$. It remains to consider the case when $G$ contains a single even cycle. On the one hand, Lemma~\ref{lemma:even} shows that $\theta_{DP}(G, k)^{-1} > 2$ for all $k \in \N^+$. On the other hand, $G$ is obtained from an even cycle by repeatedly adding vertices of degree $1$, so we can combine the result of Lemma~\ref{lemma:even} with the following obvious observation to conclude that $\chi_{DP}^\ast(G) = 2$:
		
		\begin{obs}
			If $u \in V(G)$ satisfies
			$
			\deg_G(u) \leq \chi_{DP}^\ast(G - u) - 1
			$,
			then $\chi_{DP}^\ast(G) = \chi_{DP}^\ast(G - u)$. \qedhere
		\end{obs}
	\end{scproof}
	
	\section{Proof of Theorem~\ref{theo:frac_deg}}\label{sec:frac_deg}
	

	\subsection{Notation}
	
	Here we review some notation that shall be used throughout this section. Let $G$ be a graph. For a vertex $u \in V(G)$, $N_G(u)$ and $N_G[u]$ are the open and the closed neighborhoods of $u$ respectively, i.e.,
	\[
		N_G(u) \defeq \set{v \in V(G) \,:\, uv \in E(G)} \qquad \text{and} \qquad N_G[u] \defeq N_G(u) \cup \set{u}.
	\]
	For a set $U \subseteq V(G)$, let
	\[
		N_G(U) \defeq \bigcup_{u \in U} N_G(u) \qquad \text{and} \qquad N_G[U] \defeq \bigcup_{u \in U} N_G[u] =  N_G(U) \cup U.
	\]
	Similarly, if $D$ is a digraph and $u \in V(D)$, then we write
	\begin{align*}
	N_D^+(u) &\defeq \set{v \in V(D) \,:\, uv \in E(D)};&&N_D^-(u)\defeq \set{v \in V(D) \,:\, vu \in E(D)};\\
	N_D^+[u] &\defeq N_D^+(u) \cup \set{u};&&N_D^-[u] \defeq N_D^-(u) \cup \set{u}.
	\end{align*}
	We use $R^+_D(u)$ to denote the set of all vertices $v \in V(D)$ that are {reachable} from $u$ via a directed path of positive length (so $D$ is acyclic if and only if $u \not \in R^+_D(u)$ for all $u \in V(D)$). Let $R^-_D(u)$ denote the set of all $v \in V(D)$ such that $u \in R^+_D(v)$. We write
	\[
		R^+_D[u] \defeq R^+_D(u) \cup \set{u} \qquad \text{and} \qquad R^-_D[u] \defeq R^-_D(u) \cup \set{u}.
	\]
	For a subset $U \subseteq V(D)$, let
	\[
		N_D^+(U) \defeq \bigcup_{u \in U} N_D^+(u); \qquad R_D^+(U) \defeq \bigcup_{u \in U} R_D^+(u); \qquad \text{etc}.
	\]	
	We use expressions $|S|$ and $\#S$ for the cardinality of a set $S$ interchangeably (usually, $\#S$ suggests that it is a random variable).
	
	\subsection{The random greedy algorithm}
	
	Now we can begin the proof. Let $G$, $D$, and $d$ be as in the statement of Theorem~\ref{theo:frac_deg}. For brevity, we set $V \defeq V(G)$ and omit subscripts $G$ and $D$ in expressions such as $N_G(u)$, $R^-_D[u]$, $\deg^+_D(u)$, etc. We will often use the acyclicity of $D$ to make inductive definitions or arguments by describing how to deal with a vertex $u$ provided that all $v$ reachable from $u$ have already been considered.
	
	Fix $\epsilon \in (0,1)$ and define $\eta \defeq (1-\epsilon) \ln d/d$. We will show that $\chi^\ast_{DP}(G) \leq \eta^{-1}$ if $d$ is large enough (as a function of $\epsilon$). Let $\Cov{H} = (L, H)$ be a $k$-fold cover of $G$. Our aim is to show that if $k$ is sufficiently large (where the lower bound may depend on the entire graph $G$), then $G$ has an $(\eta, \Cov{H})$-coloring.
	
	For a set $U \subseteq V$, let $L(U) \defeq \bigcup_{u \in U} L(u)$
	and let $\mathbf{QI}(U)$ denote the set of all quasi-independent sets contained in $L(U)$. Let $F$ be the orientation of the cross-edges of $H$ in which a cross-edge~$xy$ is directed from $x$ to $y$ if and only if the vertices $u$, $v \in V$ such that $x \in L(u)$ and $y \in L(v)$ satisfy $uv \in E(D)$. Again, we omit subscripts $H$ and $F$ in expressions such as $N_H[x]$, $R^+_F(x)$, etc.
	
	Given a set of probabilities $p(u) \in [0,1]$ for $u \in V$, we define random subsets $S(u) \subseteq L(u)$ inductively using the following random greedy construction. Consider $u \in V$ and suppose that the sets $S(v)$ for all $v$ reachable from $u$ have already been defined. Independently for each $x \in L(u)$, set
	\begin{equation}\label{eq:xi}
	\xi(x) \defeq \begin{cases}
	1 &\text{with probability } p(u);\\
	0 &\text{with probability } 1-p(u),
	\end{cases}
	\end{equation}
	Define
	\[
		L'(u) \defeq \set{x \in L(u) \,:\, N^+(x) \cap S(v) = \0 \text{ for all } v \in N^+(u)},
	\]
	and then
	\[
		S(u) \defeq \set{x \in L'(u) \,:\, \xi(x) = 1}.
	\]
	Note that for every $u \in V$, the set $S(u)$ only depends on the random choices associated with the elements of $L(R^+[u])$. For each $U \subseteq V$, write
	$
		S(U) \defeq \bigcup_{u \in U} S(u)
	$
	and set $S \defeq S(V)$. By construction, $S$ is always a quasi-independent set. In the remainder of this section we will analyze this construction in order to argue that, for a suitable choice of probabilities $\set{p(u)\,:\,u \in V}$ and sufficiently large $k$, $|S(u)| \geq \eta k$ for all $u \in V$ \ep{and hence $S$ is an $(\eta, \Cov{H})$-coloring} with high probability. 
	
	\subsection{A correlation inequality}
	
	Recall that a family $\mathcal{F}$ of subsets of a set $I$ is \emph{increasing} if whenever $X_1 \supseteq X_2 \in \mathcal{F}$, we have $X_1 \in \mathcal{F}$; similarly, $\mathcal{F}$ is \emph{decreasing} if $X_1 \subseteq X_2 \in \mathcal{F}$ implies $X_1 \in \mathcal{F}$. We shall need the following version of the FKG inequality, tracing back to the work of Kleitman~\cite{Kleitman}:
	
	\begin{theo}[{\cite[Theorem~6.3.2]{AS00}}]\label{theo:FKG}
		Let $I$ be a finite set and let $X$ be a random subset of $I$ obtained by selecting each $i \in I$ independently with probability $q(i) \in [0,1]$. If $\mathcal{F}$ and $\mathcal{G}$ are increasing families of subsets of~$I$, then
		\[
		\Pr\left[X \in \mathcal{F} \text{ and } X \in \mathcal{G}\right] \geq \Pr\left[X \in \mathcal{F}\right]\cdot \Pr\left[X \in \mathcal{G}\right].
		\]
		The same conclusion holds if $\mathcal{F}$ and $\mathcal{G}$ are decreasing.
	\end{theo}

	We use Theorem~\ref{theo:FKG} to obtain the following positive correlation result:
	
	\begin{lemma}\label{lemma:correlation}
		Let $u \in V$ and set
		$
			A \defeq R^+[u] \setminus R^-[N^+[u]]
		$.
		Let $Q \in \mathbf{QI}(A)$ and $Y \subseteq L(N^+(u))$. Then
		\[
			\Pr \left[y \not \in S \text{ for all } y \in Y\,\middle\vert\, S(A) = Q\right] \geq \prod_{y \,\in\, Y} \Pr\left[y \not \in S\,\middle\vert\, S(A) = Q\right].
		\]
	\end{lemma}
	\begin{scproof}\stepcounter{ForClaims} \renewcommand{\theForClaims}{\ref{lemma:correlation}}
		Since nothing in the statement of the lemma depends on the vertices outside of $R^+[u]$, we may pass to a subgraph and assume that $V = R^+[u]$. Set
		$
			B \defeq R^-[N^+[u]]
		$,
		so $B = V \setminus A$. The lemma is trivially true if $Y = \0$, so we may assume $Y \neq \0$, and hence $N^+(u) \neq \0$.
		
		The graph $G[B]$ is bipartite. Indeed, consider any $v \in B$. On the one hand, $v$ is reachable from $u$; on the other hand, there is a vertex $w \in N^+(u)$ reachable from $v$. We claim that if $P_1$ and $P_2$ are two directed $uv$-paths, then $\operatorname{length}(P_1) \equiv \operatorname{length}(P_2) \pmod 2$. Indeed, let $P_3$ be any directed $vw$-path. If $\operatorname{length}(P_1) \not\equiv \operatorname{length}(P_2) \pmod 2$, then either $P_1+P_3$ or $P_2 + P_3$ is a directed $uw$-path of even length, which contradicts assumption \ref{item:parity}. Thus, we can $2$-color the vertices in $B$ based on the parity of the directed paths leading from $u$ to them.
		
		Let $\set{U_1, U_2}$ be a partition of $B$ into two independent sets such that $u \in U_1$. To the random variable $\xi$ defined in \eqref{eq:xi}, we associate a random subset $X_\xi \subseteq L(B)$ as follows:
		\[
			X_\xi \defeq \set{x \in L(U_1)\,:\, \xi(x) = 1} \cup \set{x \in L(U_2) \,:\, \xi(x) = 0}.
		\]
		The set $X_\xi$ is built by independently selecting each element $x \in L(B)$ with probability $q(x)$ given by
		\[
			q(x) \defeq \begin{cases}
				p(v) &\text{if } x \in L(v) \text{ for } v \in U_1;\\
				1-p(v) &\text{if } x \in L(v) \text{ for } v \in U_2.
			\end{cases}
		\]
		To complete the construction of the set $S$, given that $S(A) = Q$, we only need to know the values $\xi(x)$ for all $x \in L(B)$. Since all of them are determined by the set $X_\xi$, we may, for fixed $X \subseteq L(B)$, denote by $S_X$ the value $S$ would take under the assumptions $S(A) = Q$ and $X_\xi = X$. For each $x \in L(B)$, let
		\[
			\mathcal{F}_x \defeq \set{X \subseteq L(B) \,:\, x \not\in S_X}.
		\]
		Since the random set $X_\xi$ is generated independently from $S(A)$, we have
		\[
			\Pr[X_\xi \in \mathcal{F}_x] \cdot \Pr[S(A) = Q] = \Pr[X_\xi \in \mathcal{F}_x \,\wedge\, S(A) = Q] = \Pr[x \not\in S \,\wedge\, S(A) = Q],
		\]
		and hence
		\begin{equation}\label{eq:cond_prob}
			\Pr[X_\xi \in \mathcal{F}_x] = \Pr\left[x \not\in S \,\middle\vert\, S(A) = Q\right].
		\end{equation}
		
		\begin{claim}\label{claim:monotone}
			If $x \in L(U_1)$ \ep{resp.\ $x \in L(U_2)$}, then the family $\mathcal{F}_x$ is decreasing \ep{resp.\ increasing}.
		\end{claim}
		\begin{claimproof}
			We argue inductively. Let $v \in B$ and suppose that the claim has been verified for all $x \in L(w)$ with $w \in B$ reachable from $v$. Consider any $x \in L(v)$. We will give the proof for the case $v \in U_1$, as the case $v \in U_2$ is analogous. By definition,
			\begin{equation}\label{eq:x_xi}
				x \in S \,\Longleftrightarrow\, \xi(x) = 1 \text{ and } y \not\in S \text{ for all } y \in N^+(x).
			\end{equation}
			If $x \in N^-(Q)$, then $\mathcal{F}_x = \powerset{L(B)}$ \ep{and we are done}, while if $x \not \in N^-(Q)$, then~\eqref{eq:x_xi} yields
			\[
				\overline{\mathcal{F}_x} = \set{X \subseteq L(B) \,:\, x \in X} \cap \bigcap_{y \,\in\, N^+(x) \cap L(B)} \mathcal{F}_y,
			\]
			where $\overline{\mathcal{F}_x}$ denotes the complement of $\mathcal{F}_x$. Each $y \in N^+(x) \cap L(B)$ belongs to $L(U_2)$, so, by the inductive assumption, the family $\mathcal{F}_y$ are increasing. The intersection of increasing families is again increasing, so $\mathcal{F}_x$ is the complement of an increasing family, hence $\mathcal{F}_x$ is decreasing, as desired.
		\end{claimproof}
		
		With Claim~\ref{claim:monotone} in hand, the conclusion of the lemma follows from equation~\eqref{eq:cond_prob} and the fact that $Y \subseteq L(U_2)$ by applying Theorem~\ref{theo:FKG} to the increasing families $\mathcal{F}_y$ with $y \in Y$.
	\end{scproof}

	We remark that the proof of Lemma~\ref{lemma:correlation} is the only place where we use condition \ref{item:parity}.
	
	\subsection{Expectation bounds}
	
	The next lemma gives a lower bound on the \emph{expected} sizes of the sets $S(u)$.
	
	\begin{lemma}\label{lemma:expectations}
		Let $\alpha$ be a positive real number such that
		\[
			(1+\alpha)^2(1-\epsilon) < 1.
		\]
		Then there exists a choice of $\set{p(u) \,:\, u \in V}$ such that for all $u \in V$,
		\[
			\mathbb{E}[\#S(u)] = (1+\alpha)\eta k.
		\]
	\end{lemma}
	\begin{scproof}
		Let $\beta \in (0,1)$ be such that
		\begin{equation}\label{eq:beta}
			1 - \lambda \geq \exp\left(-(1+\alpha)\lambda\right) \text{ for all } 0 < \lambda \leq \beta.
		\end{equation}
		We will frequently use the following form of the inequality of arithmetic and geometric means: Given nonnegative real numbers $\lambda_1$, \ldots, $\lambda_m$ and nonnegative weights $w_1$, \ldots, $w_m$ satisfying $\sum_{i=1}^m w_i = 1$,
		\begin{equation}\label{eq:AGM}
			\sum_{i=1}^m w_i\lambda_i \geq \prod_{i=1}^m \lambda_i^{w_i}.
		\end{equation}
		
		We define the values $p(u)$ inductively. Let $u \in V$ and assume that we have already defined $p(v)$ for all $v$ reachable from $u$ so that
		\[
			\mathbb{E}[\#S(v)] = (1+\alpha) \eta k \quad\text{and}\quad p(v) \leq \beta \quad \text{for all}\quad v \in R^+(u).
		\]
		We will show that in that case
		\begin{equation}\label{eq:good_exp}
			\mathbb{E}[\#L'(u)] \geq \beta^{-1}(1+\alpha)\eta k.
		\end{equation}
		After~\eqref{eq:good_exp} is established, we can define
		\[
			p(u) \defeq \frac{(1+\alpha)\eta k}{\mathbb{E}[\#L'(u)]},
		\]
		which gives
		\[
			\mathbb{E}[\#S(u)] = p(u) \mathbb{E}[\#L'(u)] = (1+\alpha)\eta k,
		\]
		as desired, and, furthermore, $p(u) \leq \beta$, allowing the induction to continue.
		
		As in the statement of Lemma~\ref{lemma:correlation}, let \[A \defeq R^+[u] \setminus R^-[N^+[u]].\] We have
		\begin{equation}\label{eq:mean_decomposition}
			\mathbb{E}[\# L'(u)] = \sum_{Q \,\in\, \mathbf{QI}(A)} \mathbb{E}\left[\# L'(u)\,\middle\vert\,S(A) = Q\right] \cdot \Pr\left[S(A) = Q\right].
		\end{equation}
		Consider any $Q \in \mathbf{QI}(A)$. By the linearity of expectation,
		\begin{align*}
			\mathbb{E}\left[\#L'(u)\,\middle\vert\,S(A) = Q\right] &= \sum_{x \,\in\, L(u)} \Pr\left[x \in L'(u)\,\middle\vert\, S(A) = Q\right] \\
			&=\sum_{x \,\in\, L(u)} \Pr\left[ y \not \in S \text{ for all } y \in N^+(x) \,\middle\vert\, S(A) = Q \right].
		\end{align*}
		From Lemma~\ref{lemma:correlation} we derive
		\[
			\sum_{x \,\in\, L(u)} \Pr\left[ y \not \in S \text{ for all } y \in N^+(x) \,\middle\vert\, S(A) = Q \right] \geq \sum_{x \,\in\, L(u)} \prod_{y \,\in\, N^+(x)} \Pr\left[ y \not \in S \,\middle\vert\, S(A) = Q \right],
		\]
		which, by~\eqref{eq:AGM}, is at least
		\[
			k \left(\prod_{x \,\in\, L(u)} \prod_{y \,\in\, N^+(x)} \Pr\left[ y \not \in S \,\middle\vert\, S(A) = Q \right]\right)^{1/k}.
		\]
		After changing the order of multiplication, we get
		\[
			\prod_{x \,\in\, L(u)} \prod_{y \,\in\, N^+(x)} \Pr\left[ y \not \in S \,\middle\vert\, S(A) = Q \right] \geq \prod_{v \,\in\, N^+(u)} \prod_{y \,\in\, L(v)} \Pr \left[y \not \in S\,\middle\vert\, S(A) = Q\right].
		\]
		Now consider any $v \in N^+(u)$. Let
		\[
			A_v \defeq A \cup R^+(v).
		\]
		Since $A \subseteq A_v$, the set $S(A)$ is determined by $S(A_v)$, and hence
		\[ 
			\prod_{y \,\in\, L(v)} \Pr \left[y \not \in S\,\middle\vert\, S(A) = Q\right] = \prod_{y \,\in\, L(v)} \sum_{R \,\in\, \mathbf{QI}(A_v)} \Pr\left[ y \not \in S \,\middle\vert\, S(A_v) = R \right] \cdot \Pr\left[S(A_v) = R\,\middle\vert\,S(A) = Q\right].
		\] 
		Applying~\eqref{eq:AGM} again, we see that the last expression is at least
		\begin{align}
			 &\prod_{y \,\in\, L(v)} \prod_{R \,\in\, \mathbf{QI}(A_v)} \left(\Pr\left[ y \not \in S \,\middle\vert\, S(A_v) = R \right]\right)^{\Pr\left[S(A_v) = R\,\middle\vert\,S(A) = Q\right]} \nonumber\\
			=&\prod_{R \,\in\, \mathbf{QI}(A_v)} \left(\prod_{y \,\in\, L(v)}\Pr\left[ y \not \in S \,\middle\vert\, S(A_v) = R \right] \right)^{\Pr\left[S(A_v) = R\,\middle\vert\,S(A) = Q\right]}. \label{eq:product}
		\end{align}
		Note that the set $L'(v)$ is completely determined by $S(A_v)$. This allows us to introduce notation $L'_R(v)$ for the value of $L'(v)$ under the assumption $S(A_v) = R$; or, explicitly,
		\[
			L'_R(v) \defeq L(v) \setminus N^-(R).
		\]
		Since $v \not \in A_v$, for fixed $R \in \mathbf{QI}(A_v)$ and $y \in L(v)$, we have
		\[
			\Pr\left[ y \not \in S \,\middle\vert\, S(A_v) = R \right] = \begin{cases}
				1 - p(v) &\text{if } y \in L'_R(v);\\
				1 &\text{otherwise}.
			\end{cases}
		\]
		Therefore,
		\[
			\prod_{y \,\in\, L(v)}\Pr\left[ y \not \in S \,\middle\vert\, S(A_v) = R \right] = (1 - p(v))^{|L'_R(v)|}.
		\]
		Plugging this into~\eqref{eq:product}, we obtain
		\begin{align*}
			\prod_{y \,\in\, L(v)} \Pr \left[y \not \in S\,\middle\vert\, S(A) = Q\right] &\geq \prod_{R \,\in\, \mathbf{QI}(A_v)} (1-p(v))^{|L'_R(v)| \cdot \Pr\left[S(A_v) = R\,\middle\vert\,S(A) = Q\right]}\\
			&= (1-p(v))^{\sum_{R \,\in\, \mathbf{QI}(A_v)} |L'_R(v)| \cdot \Pr\left[S(A_v) = R\,\middle\vert\,S(A) = Q\right]} \\
			&= (1-p(v))^{\mathbb{E}\left[\#L'(v)\,\middle\vert\, S(A) = Q\right]}.
		\end{align*}
		Since, by our assumption, $p(v) \leq \beta$, inequality~\eqref{eq:beta} yields
		\begin{align*}
			(1-p(v))^{\mathbb{E}\left[\#L'(v)\,\middle\vert\, S(A) = Q\right]} &\geq \exp\left(-(1+\alpha) p(v) \mathbb{E}\left[\#L'(v)\,\middle\vert\, S(A) = Q\right] \right) \\
			&= \exp\left(-(1+\alpha) \mathbb{E}\left[\#S(v)\,\middle\vert\, S(A) = Q\right] \right).
		\end{align*}
		This allows us to lower bound $\mathbb{E}\left[\#L'(u)\,\middle\vert\,S(A) = Q\right]$ as
		\begin{align*}
			\mathbb{E}\left[\#L'(u)\,\middle\vert\,S(A) = Q\right] &\geq k \left(\prod_{v \,\in\, N^+(u)} \exp\left(-(1+\alpha) \mathbb{E}\left[\#S(v)\,\middle\vert\, S(A) = Q\right] \right)\right)^{1/k}\\
			&= k \exp\left(-\frac{1+\alpha}{k} \sum_{v \,\in\, N^+(u)} \mathbb{E}\left[\#S(v)\,\middle\vert\, S(A) = Q\right]\right).
		\end{align*}
		Returning to~\eqref{eq:mean_decomposition}, we conclude
		\[
			\mathbb{E}[\# L'(u)] \geq k \sum_{Q \,\in\, \mathbf{QI}(A)} \exp\left(-\frac{1+\alpha}{k} \sum_{v \,\in\, N^+(u)} \mathbb{E}\left[\#S(v)\,\middle\vert\, S(A) = Q\right]\right) \cdot \Pr\left[S(A) = Q\right].
		\]
		Due to the convexity of the exponential function (or by~\eqref{eq:AGM} again), the last expression is at least
		\[
			k \exp\left(-\frac{1+\alpha}{k} \sum_{v \,\in\, N^+(u)} \mathbb{E}\left[\#S(v)\right]\right),
		\]  
		which, since $\deg^+(u) \leq d$ and $\mathbb{E}[\#S(v)] = (1+\alpha) \eta k$ for all $v \in N^+(u)$ by assumption, finally becomes
		\[
			k\exp\left(-(1+\alpha)^2\eta \deg^+(u)\right) \geq k\exp\left(-(1+\alpha)^2\eta d\right) = k d^{-(1+\alpha)^2(1-\epsilon)}.
		\]
		It remains to notice that, since $(1+\alpha)^2(1-\epsilon) < 1$, the quantity $d^{-(1+\alpha)^2(1-\epsilon)}$ is asymptotically larger than $\beta^{-1}(1+\alpha)\eta = \Theta(\ln d/d)$. This finishes the proof of~\eqref{eq:good_exp}.
	\end{scproof}

	\subsection{Concentration of measure and completing the proof}
	
	Finally, we show that the sizes of the sets $S(u)$ are highly concentrated around their expected values. The tool that we use is the following classical result, which is a consequence of
	Azuma's inequality for Doob martingales \ep{see \cite[\S7.4]{AS00}}:
	
	\begin{theo}[{\cite[79]{MR02}}]\label{theo:SCB}
		Let $\zeta$ be a random variable determined by $s$ independent trials such that changing the outcome of any one trial can affect the value of $\zeta$ at most by~$c$. Then 
		\[
		\Pr[|\zeta - \mathbb{E} \zeta| > t] \leq 2\exp\left(-\frac{t^2}{2c^2s}\right).
		\]
	\end{theo}

	\begin{lemma}\label{lemma:concentration}
		There is $C > 0$, depending on $G$ but not on $k$, such that for all $\alpha > 0$ and $u \in V$,
		\[
			\Pr\left[\left|\#S(u) - \mathbb{E}[\#S(u)]\right| > \alpha k \right] \leq 2\exp\left(-C\alpha^2k\right).
		\]
	\end{lemma}
	\begin{scproof}\stepcounter{ForClaims} \renewcommand{\theForClaims}{\ref{lemma:concentration}}
		The value $\#S(u)$ is determined by $k|V|$ independent trials, namely by the values $\xi(x)$ for $x \in V(H)$, so, to apply Theorem~\ref{theo:SCB}, we only need to establish the following:
		
		\begin{claim}\label{claim:Lipschitz}
			Changing the value $\xi(x)$ for some $x \in V(H)$ can affect $\#S(u)$ at most by some amount~$c$ that depends on $G$ but not on $k$.
		\end{claim}
		\begin{claimproof}
			Suppose that $x \in L(v)$ for some $v \in V$. The value $\xi(x)$ can only affect $y \in R^-[x]$, so it suffices to bound $|R^-[x] \cap L(u)|$ from above. Let $y = z_1 \to \cdots \to z_\ell = x$ be a directed $yx$-path for some $y \in L(u)$. For each $1 \leq i \leq \ell$, choose $v_i \in V$ so that $z_i \in L(v_i)$. Then $u = v_1 \to \cdots \to v_\ell = v$ is a directed $uv$-path in~$D$. Notice that the $uv$-path $v_1 \to \cdots \to v_\ell$ uniquely identifies $y = z_1$. Indeed, by definition, $z_\ell = x$, so $z_{\ell-1}$ must be the unique neighbor of $x$ in $L(v_{\ell-1})$. Then $z_{\ell-2}$ must be the unique neighbor of $z_{\ell-1}$ in $L(v_{\ell-2})$; and so on. Thus, $|R^-[x] \cap L(u)|$ does not exceed the number of directed $uv$-paths, which is independent of $k$.
		\end{claimproof}
		
		The conclusion of the lemma now follows by applying Theorem~\ref{theo:SCB} with $s = k |V|$, $t = \alpha k $, and the value of $c$ given by Claim~\ref{claim:Lipschitz}.
	\end{scproof}
	
	Now we can easily finish the proof of Theorem~\ref{theo:frac_deg}. Pick some $\alpha > 0$ so that $(1+\alpha)^2(1-\epsilon) < 1$ and apply Lemma~\ref{lemma:expectations} to obtain $\set{p(u) \,:\, u \in V}$ such that for all $u \in V$,
	\[
		\mathbb{E}[\#S(u)] = (1+\alpha)\eta k.
	\]
	Then, by Lemma~\ref{lemma:concentration}, we have
	\[
		\Pr\left[\#S(u) < \eta k\right] \leq \left[\left|\#S(u) - \mathbb{E}[\#S(u)]\right| > \alpha\eta k \right] \leq 2\exp\left(-C\alpha^2\eta^2k\right),
	\]
	where $C > 0$ may depend on $G$ but not on $k$. Therefore, we may apply the union bound to get
	\[
		\Pr\left[\text{$S$ is not an $(\eta, \Cov{H})$-coloring of $G$}\right] \leq 2n\exp\left(-C\alpha^2\eta^2k\right) \xrightarrow[k \to \infty]{} 0,
	\]
	as desired.
	
	\section{Proof of Corollary~\ref{corl:large_girth}}\label{sec:BD}
	
	For $t \in \N^+$, a \emph{$t$-uniform hypergraph} $H$ is a pair $(V(H), E(H))$, where $V(H)$ is a set whose elements are called the \emph{vertices} of $H$ and $E(H)$ is a set of $t$-element subsets of $V(H)$, called the \emph{edges} of $H$. \ep{Thus, a graph is a $2$-uniform hypergraph.} A hypergraph $H$ is \emph{$k$-colorable} if there is a function $f \colon V(H) \to [k]$ such that for every edge $e \in E(H)$, the restriction of $f$ onto $e$ is not constant. A \emph{\ep{Berge} cycle} of length $g$ in a hypergraph $H$ consists of a sequence of distinct vertices $v_1$, \ldots, $v_g$ and a sequence of distinct edges $e_1$, \ldots, $e_g$ with the property that $e_i \supseteq \set{v_i, v_{i+1}}$ for all $1 \leq i \leq g$, where the indices are taken modulo $g$. The \emph{girth} of a hypergraph $H$ is the least integer $g \geq 2$ such that $H$ contains a cycle of length $g$ \ep{if no such $g$ exists, then the girth of $H$ is set to be $\infty$}. The following fact is well-known:
	
	\begin{theo}[{Erd\H{o}s--Hajnal \cite{EH}}]\label{theo:EH}
		For all $t$, $k$, $g \in \N^+$, there is a finite non-$k$-colorable $t$-uniform hypergraph of girth at least $g$.
	\end{theo}
	
	In \cite{BD}, Blanche Descartes introduced a simple way of building graphs with girth $6$ and arbitrarily high chromatic number. This construction was generalized in \cite{KN99} using Theorem~\ref{theo:EH} to obtain graphs whose girth and chromatic number are both arbitrarily high. The graphs produced by this construction also serve as examples for Corollary~\ref{corl:large_girth}. Start by setting $G_1 \defeq K_2$ and let $D_1$ be an orientation of $G_1$. Once $G_i$ and $D_i$ are defined, take a $|V(G_i)|$-uniform non-$(i+1)$-colorable hypergraph $H_i$. Build $G_{i+1}$ by making $V(H_i)$ an independent set, adding $|E(H_i)|$ disjoint copies of $G_i$, establishing a bijection between the copies of $G_i$ and the edges of $H_i$, and joining each copy to its corresponding edge via a perfect matching. Finally, let $D_{i+1}$ be the orientation of $G_{i+1}$ obtained by orienting each copy of $G_i$ according to $D_i$ and directing every remaining edge toward its endpoint in $V(H_i)$.
	
	The graphs $G_i$ have the following properties \ep{see \cite{KN99}}:
	
	
	\begin{itemize}
		\item[--] $G_i$ is $i$-degenerate;
		
		
		\item[--] $\chi(G_i) = i + 1$;
		
		
		\item[--] if for all $j < i$, the girth of $H_j$ is at least $g$, then the girth of $G_i$ is at least $3g$. 
	\end{itemize}
	
	Additionally, it is clear from the construction that the orientation $D_i$ is acyclic and the out-degree of every vertex in $D_i$ is at most $i$; in other words, $D_i$ witnesses that $G_i$ is $i$-degenerate. Furthermore, the \ep{undirected} subgraph of $G_i$ induced by the vertices reachable in $D_i$ from any given vertex $u \in V(G_i)$, including $u$ itself, is acyclic; in particular, for all $uv \in E(D_i)$, the only directed $uv$-path is the single edge $u \to v$. Therefore, condition~\ref{item:parity} of Theorem~\ref{theo:frac_deg} holds and we can conclude $\chi^\ast_{DP}(G_i) \leq (1+o(1))i / \ln i$. Hence, by using hypergraphs $H_i$ of large girth in this construction, we obtain graphs satisfying all the requirements of Corollary~\ref{corl:large_girth}.
	
	\subsubsection*{Acknowledgements}
	
	The authors are grateful to Doug West for helpful discussions and to the anonymous referees for their comments and suggestions.
	
	{\renewcommand{\markboth}[2]{}
		\printbibliography}
	
\end{document}